\documentclass[11pt]{article}
\usepackage{amssymb,amsfonts,amsmath,latexsym, epsfig,mathrsfs,colordvi}
\usepackage{graphicx}
\newtheorem{theo}{Theorem}

\newtheorem{lem} [theo]{Lemma}
\newtheorem{coro}[theo]{Corollary}

\makeatletter \@addtoreset{equation}{section}
\@addtoreset{theo}{}\makeatother

\def\qed{\hfill \rule{4pt}{7pt}}
\def\pf{\noindent {\it Proof.} }

\begin{document}

\title{Counting Humps in Motzkin paths}

\author{Yun Ding and Rosena R. X. Du\footnote{Email: rxdu@math.ecnu.edu.cn.} \\ \\ Department of Mathematics, East China Normal University\\
500 Dongchuan Road, Shanghai, 200241, P. R. China.}

\date{August 20, 2011}
\maketitle

\vskip 0.7cm \noindent{\bf Abstract.}
In this paper we study the number of humps (peaks) in Dyck, Motzkin and Schr\"{o}der paths.
Recently A. Regev noticed that the number of peaks in all Dyck paths of order $n$ is one half of the number of super Dyck paths of order $n$. He also computed the number of humps in Motzkin paths and found a similar relation, and asked for bijective proofs. We give a bijection and prove these results. Using this bijection we also give a new proof that the number of Dyck paths of order $n$ with $k$ peaks is the Narayana
number. By double counting super Schr\"{o}der
paths, we also get an identity involving products of
binomial coefficients.



\section{Introduction}

A \emph{Dyck path} of order (semilength) $n$ is a lattice path in
$\mathbb{Z}\times\mathbb{Z}$, from $(0,0)$ to $(2n,0)$, using
up-steps $(1,1)$ (denoted by $U$) and down-steps $(1,-1)$ (denoted by
$D$) and never going below the $x$-axis. We use $\mathcal {D}_{n}$ to denote the set of Dyck paths of order $n$. It is well known that $\mathcal {D}_{n}$ is counted by the $n$-th \emph{Catalan number} (A000108 in \cite{sequence})
\[C_{n}=\frac{1}{n+1}{{2n}\choose{n}}.\]

A \emph{peak} in a Dyck path is two consecutive steps $UD$. It is also well known (see, for example, \cite{Deutsch, Narayana, stanleyec2}) that the number of Dyck paths of order $n$ with $k$
peaks is the \emph{Narayana number} (A001263):
\[N(n,k)=\frac{1}{n}{{n}\choose{k}}{{n}\choose{k-1}}.\]
Counting Dyck paths with restriction on peaks has been studied by many authors, see for example \cite{Mansour1, Mansour2, PeartWoan}. Here we are interested in counting peaks in all Dyck paths of order $n$. By summing over the above formula over $k$ we immediately get the following result: the total number of peaks in all Dyck paths of order $n$ is
\[pd_{n}=\sum_{k=1}^{n}kN(n,k)={{2n-1}\choose{n}}.\]

If we allow a Dyck path to go bellow the $x$-axis, we get a \emph{super Dyck path}. Let $\mathcal{SD}_{n}$ denote the set of super Dyck paths of order $n$. By standard arguments we have
\begin{equation}\label{2Dyck}
sd_{n}=\#\mathcal{SD}_{n}={{2n}\choose{n}}=2{{2n-1}\choose{n}}=2 pd_{n},
\end{equation}
That is, the number of super Dyck paths of order $n$ is twice the number of peaks in all Dyck paths of order $n$. This curious relation was first noticed by Regev \cite{HumpsRegev}, who also noticed that similar relation holds for Motzkin paths, which we will explain next.

A \emph{Motzkin path of order $n$} is a lattice path in
$\mathbb{Z}\times\mathbb{Z}$, from $(0,0)$ to $(n,0)$, using
up-steps $(1,1)$, down-steps $(1,-1)$ and flat-steps $(1,0)$ (denoted by $F$) that never goes below
the $x$-axis. Let $\mathcal {M}_{n}$ denote all
the Motzkin paths of order $n$. The cardinality of
$\mathcal {M}_{n}$ is the $n$-th \emph{Motzkin number} $m_{n}$ (A001006), which satisfies the following recurrence relation
\[m_{0}=1,~~m_{1}=1,~~m_{n}=m_{n-1}+\sum_{i=2}^{n}m_{i-2}m_{n-i}, ~~for~~ n\geq2,\]
and have generating function
\[\sum_{n\geq0}m_{n}x^{n}=\frac{1-x-\sqrt{1-2x-3x^{2}}}{2x^{2}}.\]

A \emph{hump} in a Motzkin path is an up step followed by zero or more flat steps followed by a down step. We use $hm_{n}$ to denote the total number of humps in all Motzkin paths of order $n$. We can similarly define \emph{super Motzkin paths} to be Motzkin paths that are allowed to go below the $x$-axis, and use $\mathcal{SM}_n$ to denote the set of super Motzkin paths of order $n$. Using a recurrence relation and the WZ method \cite{A=B,Zeilberger}, Regev (\cite{HumpsRegev}) proved that
\begin{equation}\label{2Motzkin}
sm_{n}=\#\mathcal{SM}_{n}=\sum_{j\geq0}{{n}\choose{j}}{{n-j}\choose{j}}=2 hm_{n}+1
\end{equation}
and asked for a bijective proof of \eqref{2Dyck} and \eqref{2Motzkin}. The main result of this paper is such a bijective proof.

Let $\mathcal{SM}_{n}^{UU}(k)$ ($\mathcal{SM}_{n}^{UD}(k)$) denote the set of paths in $\mathcal {SM}_{n}$ with $k$ peaks and the first non-flat step is $U$, and the last
non-flat step is $U$ ($D$). Let $\mathcal {SM}_{n}^{U*}$ denote all paths in $\mathcal {SM}_{n}$ whose first non-flat step is $U$, and define
\[\mathcal{HM}_{n}=\{(M,P)|M\in\mathcal{M}_{n},\text{$P$ is a hump of $M$}\}.
\]

The main result of this paper is the following:

\begin{theo}\label{th:main}
There is a bijection $\Phi: \mathcal{HM}_{n} \rightarrow \mathcal{SM}_{n}^{U*}$ such that if
$(M,P)\in\mathcal{HM}_{n}$ and $L=\Phi(M,P)$, then there are $k$ humps in $M$ if and only if $L \in\mathcal{SM}_{n}^{UU}(k-1)\cup\mathcal{SM}_{n}^{UD}(k).$
\end{theo}

The outline of the paper is as follows. In Section 2 we define the bijection $\Phi$ and prove Theorem \ref{th:main}. In section 3 we apply $\Phi$ to Dyck paths and give a new proof of the Narayana numbers. In section 4 we apply $\Phi$ to Schr\"{o}der paths and get an identity involving products of binomial coefficients by double counting super Schr\"{o}der paths whose $F$ steps are $m$-colored.

\section{The bijection $\Phi:\mathcal{HM}_{n}\leftrightarrow\mathcal{SM}_{n}^{U*}$}

Note that a Motzkin path $M$ of order $n$ can also be considered as a sequence
$M=M_{1}M_{2}\cdots M_{n},$ with $M_{i}\in\{U,F,D\},$ and the
number of $U$'s is not less than the number of $D$'s in every subsequence
$M_{1}M_{2}\cdots M_{k}$ of $M$. Hence a hump in $M$ is a subsequence $P=M_{i}M_{i+1}\cdots M_{i+k+1}, k\geq 0$, such that $M_{i}=U$, $M_{i+1}=M_{i+2}=\cdots=M_{i+k}=F$ and $M_{i+k+1}=D$. We call the end point of step $M_{i}$ a \emph{hump point}, and will also denoted as $P$. Similarly, if there exists
$i$ such that $M_{i}=D$, $M_{i+1}=M_{i+2}=\cdots=M_{i+k}=F, k\geq0$, $M_{i+k+1}=U$, then we
call the subsequence $M_{i}M_{i+1}\cdots M_{i+k+1}$ a \emph{valley} of $M$, and the end point of $M_{i+k}$ is called a
\emph{valley point}. The end point $(n,0)$ of $M$ is also considered as a
valley point.

Suppose $L$ is a path in $\mathbb{Z}\times\mathbb{Z}$ from $O(0,0)$ to $N(n,0)$, and $A$ a lattice point on $M$, we use $x_{A}$ and $y_{A}$ to denote the $x$-coordinate and $y$-coordinate of $A$, respectively. The sub-path of $L$ from point $A$ to point $B$ is denoted by
$L_{AB}$. We use $\bar{L}$ to denote the lattice path obtained from $L$ by interchanging all the up-steps and down-steps in $L$, and keep the
flat-steps unchanged.

Now we are ready to define the map $\Phi$ and prove Theorem \ref{th:main}.

\noindent {\it Proof of Theorem \ref{th:main}}:

(1) The map $\Phi:\mathcal{HM}_{n}\rightarrow\mathcal{SM}_{n}^{U*}$.

For any $(M,P)\in\mathcal{HM}_{n}$, we define $L=\Phi(M,P)$ by the following rules:
\begin{itemize}
\item  Let $C$ be the leftmost valley point in $M$ such that $x_{C}>x_{P}$;
\item  Let $B$ be the rightmost point in $M$ such that $x_{B}<x_{P},y_{B}=y_{C}$;
\item  Let $A$ be the rightmost point in $M$ such that $y_{A}=0, x_{A}\leq x_{B}$;
\item  Set $L_{0}=M_{OA}$, $L_{1}=M_{AB}$, $L_{2}=M_{BC}$, $L_{3}=M_{CN}$ (Note that $L_{0}$, $L_{1}$ and $L_{3}$ may be empty);
\item  Define $L=\Phi (M,P)=L_{0} L_{2} \overline{L_{3}} \overline{L_{1}}$.
\end{itemize}

Now we will prove that ~$L\in \mathcal {S}\mathcal {M}_{n}^{U*}$. According to the above definition, $L_0$ and $L_{2}$ are both Motzkin paths, therefore $\#U=\#D$ in $L_{0}$ and $L_{2}$. And for $L_{1}$, we have $\#U-\#D=y_{B}-y_{A}=y_{B}=y_{C}$, for $L_{3}$, $\#U-\#D=-y_{C}$. Therefore the total number of $U$'s is as much as that of  $D$'s in $L$. Thus $L$ is a super Motzkin
path of order $n$. Moreover, the first non-flat step in $L$ must be
in $L_{0}$ (when $L_{0}$ is not empty) or in $L_{2}$ (when $L_{0}$ is empty), and $L_{0}, L_{2}$ are both Motzkin paths, hence the first step leaving the $x$-axis must be a $U$. Therefore we proved that  $L=\Phi (M,P)\in\mathcal{SM}_{n}^{U*}$.

(2) The inverse of $\Phi$.

For any $L\in \mathcal {S}\mathcal {M}_{n}^{U*}$, we define $\Psi$ by the following rules:
\begin{itemize}
\item  Let $B$ be the leftmost point such that $y_{B}=0$, and $L$ goes below the $x$-axis after $B$. (If such a point does not exist, then set $B=N$);
\item  Let $A$ be the rightmost point in $L$ such that $x_{A}<x_{B}, y_{A}=0$;
\item  Let $C$ be the rightmost point in $L$ such that $x_{C}\geq x_{B}$, and $\forall G$, $x_{G}\geq x_{B}$ implies that $y_{C}\geq y_{G}$;
\item  Let $P$ be the rightmost hump point in $L$ such that $x_{P}<x_{B}$;
\item  Set $L_{0}=L_{OA}$, $L_{1}=L_{AB}$, $L_{2}=L_{BC}$, $L_{3}=L_{CN}$ (Note that $L_{0}$, $L_{2}$ and $L_{3}$ may be empty);
\item  Set $M=L_{0}\overline{L_{3}}L_{1}\overline{L_{2}}$, and $\Psi(L)=(M,P)$.
\end{itemize}

Now we prove that $\Psi=\Phi^{-1}$. Since $C$ is the highest point in $L_{3}$, and
$\overline{L_{3}}$ and $L_{3}$ are symmetric with respect to the line
$y=y_{C}$, $C$ is mapped to the lowest point in
$\overline{L_{3}}$. Moreover, $L_{0}$ and $L_{1}$ are both Motzkin
paths, then ${L_{0}}\overline{L_{3}}L_{1}$ does not go below
the $x$-axis, and the $y$-coordinate of the end point of ${L_{0}}\overline{L_{3}}L_{1}$ is
$y_{C}$. In $\overline{L_{2}}$, the end point is the lowest point, and the start point of $\overline{L_{2}}$ is $y_{C}$ higher than the end point. So
$M=L_{0}\overline{L_{3}}L_{1}\overline{L_{2}}$ ends on the $x$-axis and never goes below it, i.e., $M \in \mathcal{M}_n$. Thus $\Psi(L) \in \mathcal{HM}_n$, and it is not hard to see that $\Psi=\Phi^{-1}$.

(3) There are $k$ humps in $M$ if and only if $\Phi(M,P)\in\mathcal {SM}_{n}^{UD}(k)\cup\mathcal{SM}_{n}^{UU}(k-1).$

Since $\Phi (M)=L_{0} L_{2} \overline{L_{3}} \overline{L_{1}}=L$, the number of humps changes only in sub-paths $\overline{L_{3}}$ and $ \overline{L_{1}}$ when $M$ is converted to $L$. If the last step of $L_{1}$ is $U$, then the last step in $\overline{L_{1}}$ becomes $D$. The
number of humps in $L_{1}$ is the same as the number of humps in $\overline{L_{1}}$, and the number of
humps in $\overline{L_{3}}$ is $1$ less than the number of humps in $L_{3}$. The
last step in $\overline{L_{3}}$ is $U$ step, so concatenating
$\overline{L_{1}}$ with $\overline{L_{3}}$ yields a new hump. Therefore the total number of humps in $L$ is the same as in $M$. Thus we have $\Phi(M,P)\in\mathcal
{S}\mathcal {M}_{n}^{UD}(k).$

If the last step in $L_{1}$ is $D$, then the last step in
$\overline{L_{1}}$ is $U$. The number of humps in $\overline L_{1}$
is $1$ less than the number of humps in ${L_{1}}$, and the humps in $\overline{L_{3}}$
is $1$ less than the number of humps in $L_{3}$. Moreover, the last step in
$\overline{L_{3}}$ is $U$, so concatenating $\overline{L_{1}}$ with
$\overline{L_{3}}$ yields a new hump. Therefore the total number of humps in $L$ is $1$ less than the number humps in $M$. Thus we have $\Phi(M,P)\in\mathcal {S}\mathcal
{M}_{n}^{UU}(k-1).$

\qed

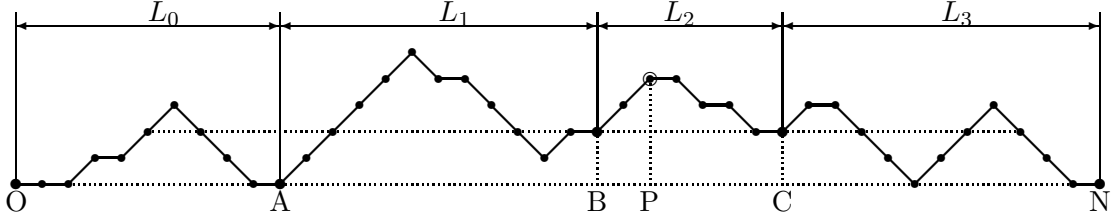
\begin{figure}[!htbp]
\centering
\begin{center}
\begin{picture}(440,70)
\put(0,0){\circle*{4}} \put(10,0){\circle*{3}}
\put(20,0){\circle*{3}} \put(30,10){\circle*{3}}
\put(40,10){\circle*{3}} \put(50,20){\circle*{3}}
\put(60,30){\circle*{3}} \put(70,20){\circle*{3}}
\put(80,10){\circle*{3}} \put(90,0){\circle*{3}}
\put(100,0){\circle*{4}} \put(110,10){\circle*{3}}
\put(120,20){\circle*{3}} \put(130,30){\circle*{3}}
\put(140,40){\circle*{3}} \put(150,50){\circle*{3}}
\put(160,40){\circle*{3}} \put(170,40){\circle*{3}}
\put(180,30){\circle*{3}} \put(190,20){\circle*{3}}
\put(200,10){\circle*{3}} \put(210,20){\circle*{3}}
\put(220,20){\circle*{4}} \put(230,30){\circle*{3}}
\put(240,40){\circle*{3}} \put(240,40){\circle{5}}
\put(250,40){\circle*{3}} \put(260,30){\circle*{3}}
\put(270,30){\circle*{3}} \put(280,20){\circle*{3}}
\put(290,20){\circle*{4}}
\put(300,30){\circle*{3}}\put(310,30){\circle*{3}}
\put(320,20){\circle*{3}} \put(330,10){\circle*{3}}
\put(340,0){\circle*{3}}\put(350,10){\circle*{3}}
\put(360,20){\circle*{3}}\put(370,30){\circle*{3}}
\put(380,20){\circle*{3}}\put(390,10){\circle*{3}}
\put(400,0){\circle*{3}}\put(410,0){\circle*{4}}

{\thinlines \put(0,0){\line(0,1){65}} \put(100,0){\line(0,1){65}}
\put(220,20){\line(0,1){45}} \put(290,20){\line(0,1){45}}
\put(410,0){\line(0,1){65}}}

\put(50,60){\vector(1,0){50}} \put(160,60){\vector(1,0){60}}
\put(350,60){\vector(1,0){60}} \put(50,60){\vector(-1,0){50}}
\put(160,60){\vector(-1,0){60}} \put(350,60){\vector(-1,0){60}}
\put(230,60){\vector(1,0){60}} \put(230,60){\vector(-1,0){10}}

\put(50,60){\makebox(0,0)[bl]{$L_{0}$}}
\put(160,60){\makebox(0,0)[bl]{$L_{1}$}}
\put(245,60){\makebox(0,0)[bl]{$L_{2}$}}
\put(350,60){\makebox(0,0)[bl]{$L_{3}$}}

\qbezier[205](0,0)(205,0)(410,0)
\qbezier[165](50,20)(215,20)(380,20)
\qbezier[10](220,0)(220,10)(220,20)
\qbezier[20](240,0)(240,20)(240,40)\qbezier[10](290,0)(290,10)(290,20)

\thicklines \put(0,0){\line(1,0){20}}\put(20,0){\line(1,1){10}}
\put(30,10){\line(1,0){10}}\put(40,10){\line(1,1){20}}
\put(60,30){\line(1,-1){30}} \put(90,0){\line(1,0){10}}
\put(100,0){\line(1,1){50}} \put(150,50){\line(1,-1){10}}
\put(160,40){\line(1,0){10}} \put(170,40){\line(1,-1){30}}
\put(200,10){\line(1,1){10}} \put(210,20){\line(1,0){10}}
\put(220,20){\line(1,1){20}} \put(240,40){\line(1,0){10}}
\put(250,40){\line(1,-1){10}} \put(260,30){\line(1,0){10}}
\put(270,30){\line(1,-1){10}} \put(280,20){\line(1,0){10}}
\put(290,20){\line(1,1){10}} \put(300,30){\line(1,0){10}}
\put(310,30){\line(1,-1){30}} \put(340,0){\line(1,1){30}}
\put(370,30){\line(1,-1){30}} \put(400,0){\line(1,0){10}}
 \put(-4,-10){\makebox(0,0)[bl]{O}}
\put(96,-10){\makebox(0,0)[bl]{A}}
\put(406,-10){\makebox(0,0)[bl]{N}}
\put(216,-10){\makebox(0,0)[bl]{B}}
\put(236,-10){\makebox(0,0)[bl]{P}}
\put(286,-10){\makebox(0,0)[bl]{C}}
\end{picture}
\end{center}
\caption{A Motzkin path $M \in \mathcal{M}_{41}$ with a circled hump point $P$. \label{Fig1}}
\end{figure}

\begin{figure}[!htbp]
\centering
\begin{center}
\begin{picture}(440,60)
\put(0,0){\circle*{4}} \put(10,0){\circle*{3}}
\put(20,0){\circle*{3}} \put(30,10){\circle*{3}}
\put(40,10){\circle*{3}} \put(50,20){\circle*{3}}
\put(60,30){\circle*{3}} \put(70,20){\circle*{3}}
\put(80,10){\circle*{3}} \put(90,0){\circle*{3}}
\put(100,0){\circle*{4}} \put(110,10){\circle*{3}}
\put(120,20){\circle*{3}} \put(130,20){\circle*{3}}
\put(140,10){\circle*{3}} \put(150,10){\circle*{3}}
\put(160,0){\circle*{3}} \put(170,0){\circle*{4}}
\put(180,-10){\circle*{3}} \put(190,-10){\circle*{3}}
\put(200,0){\circle*{3}} \put(210,10){\circle*{3}}
\put(220,20){\circle*{3}} \put(230,10){\circle*{3}}
\put(240,0){\circle*{3}} \put(250,-10){\circle*{3}}
\put(260,0){\circle*{3}} \put(270,10){\circle*{3}}
\put(280,20){\circle*{3}} \put(290,20){\circle*{4}}
\put(300,10){\circle*{3}}\put(310,0){\circle*{3}}
\put(320,-10){\circle*{3}} \put(330,-20){\circle*{3}}
\put(340,-30){\circle*{3}}\put(350,-20){\circle*{3}}
\put(360,-20){\circle*{3}}\put(370,-10){\circle*{3}}
\put(380,0){\circle*{3}}\put(390,10){\circle*{3}}
\put(400,0){\circle*{3}}\put(410,0){\circle*{4}}

\qbezier[205](0,0)(205,0)(410,0) \qbezier[10](120,0)(120,10)(120,20)
\qbezier[10](290,0)(290,10)(290,20)

{\thinlines \put(0,0){\line(0,1){45}} \put(100,0){\line(0,1){45}}
\put(170,0){\line(0,1){45}} \put(290,20){\line(0,1){25}}
\put(410,0){\line(0,1){45}}}

\put(50,40){\vector(1,0){50}} \put(130,40){\vector(1,0){40}}
\put(350,40){\vector(1,0){60}} \put(50,40){\vector(-1,0){50}}
\put(130,40){\vector(-1,0){30}} \put(350,40){\vector(-1,0){60}}
\put(200,40){\vector(1,0){90}} \put(200,40){\vector(-1,0){30}}

\put(50,40){\makebox(0,0)[bl]{$L_{0}$}}
\put(135,40){\makebox(0,0)[bl]{$L_{2}$}}
\put(225,40){\makebox(0,0)[bl]{$\overline{L_{3}}$}}
\put(350,40){\makebox(0,0)[bl]{$\overline{L_{1}}$}}

\thicklines

\put(0,0){\line(1,0){20}}\put(20,0){\line(1,1){10}}
\put(30,10){\line(1,0){10}}\put(40,10){\line(1,1){20}}
\put(60,30){\line(1,-1){30}} \put(90,0){\line(1,0){10}}
\put(100,0){\line(1,1){20}} \put(120,20){\line(1,0){10}}
\put(130,20){\line(1,-1){10}} \put(140,10){\line(1,0){10}}
\put(150,10){\line(1,-1){10}} \put(160,0){\line(1,0){10}}
\put(170,0){\line(1,-1){10}} \put(180,-10){\line(1,0){10}}
\put(190,-10){\line(1,1){30}} \put(220,20){\line(1,-1){30}}
\put(250,-10){\line(1,1){30}} \put(280,20){\line(1,0){10}}
\put(290,20){\line(1,-1){50}}\put(340,-30){\line(1,1){10}}
\put(350,-20){\line(1,0){10}} \put(360,-20){\line(1,1){30}}
\put(390,10){\line(1,-1){10}} \put(400,0){\line(1,0){10}}

\put(-4,-10){\makebox(0,0)[bl]{O}}
\put(96,-10){\makebox(0,0)[bl]{A}}
\put(406,-10){\makebox(0,0)[bl]{N}}
\put(166,-10){\makebox(0,0)[bl]{B}}
\put(116,-10){\makebox(0,0)[bl]{P}}
\put(286,-10){\makebox(0,0)[bl]{C}}
\end{picture}
\end{center}
\caption{A super Motzkin path $L=\Phi(M,P)$. \label{Fig2}}
\end{figure}
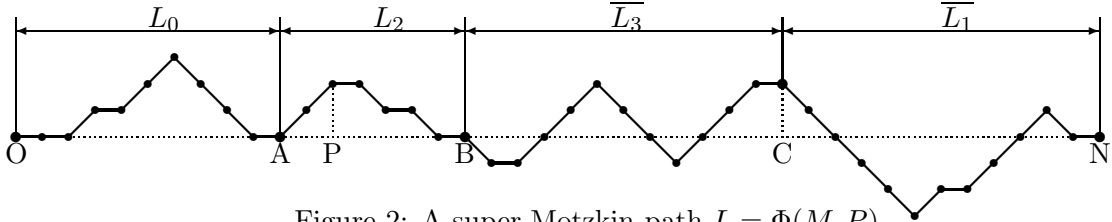

As an example, Figure \ref{Fig1} shows a Motzkin path $M \in \mathcal{M}_{41}$ with a circled hump point $P$, and Figure \ref{Fig2} shows a super Motzkin path $L \in \mathcal{SM}_{41}^{U*}=\Phi(M,P)$.

From Theorem \ref{th:main} we can easily get the following result.

\begin{coro}\label{th:theorem1} For all $n \geq 0$, we have
\begin{equation}\label{eq:smhm}
sm_{n}=2hm_{n}+1,
\end{equation}
and
\begin{equation}\label{eq:hm}
hm_{n}=\frac{1}{2}\left(\sum_{j\geq0}{{n}\choose{j}}{{n-j}\choose{j}}-1\right).
\end{equation}
\end{coro}

\pf Equation \eqref{eq:smhm} follows immediately from Theorem \ref{th:main}. To prove \eqref{eq:hm} we count super Motzkin paths with $j$ $U$ steps.
We can first choose the $j$ $U$ steps among the total $n$ steps, then
choose $j$ steps as $D$ steps among the remaining $n-j$ steps. Thus we have
\[
sm_{n}=\sum_{j\geq0}{{n}\choose{j}}{{n-j}\choose{j}}.
\]
Combine with equation \eqref{eq:smhm} we get equation
\eqref{eq:hm}. \qed

\section{Counting peaks in Dyck paths and the Narayana numbers}

Note that when restricted to Dyck paths, $\Phi$ is a bijection between super Dyck paths and peaks in Dyck paths. Therefore we have the following result.

\begin{coro}\label{th:pdn} For all $n \geq 0$, we have
\begin{equation*}
sd_{n}=2pd_{n},
\end{equation*}
and
\begin{equation*}
pd_{n}={{2n-1}\choose{n}}.
\end{equation*}
\end{coro}

Moreover, from the bijection $\Phi$ we can easily get a new proof for the Narayana numbers. To this end we need the following lemma.

\begin{lem}\label{yinli}
Let $\mathcal {SD}_{n}^{UD}(k)$ ($\mathcal {SD}_{n}^{UU}(k)$) denote the set of super Dyck paths of order $n$ with $k$ peaks whose first step is $U$ and last step is $D$ ($U$), then we have
\begin{eqnarray}
\#\mathcal {SD}_{n}^{UD}(k)={{n-1}\choose{k-1}}^{2},\label{udpeak}\\
\#\mathcal{SD}_{n}^{UU}(k)={{n-1}\choose{k-1}}{{n-1}\choose{k}}\label{uupeak},
\end{eqnarray}
and the number of super Dyck paths with $k$ peaks of order $n$ is
${{n}\choose{k}}^{2}.$
\end{lem}

\pf Each $L\in\mathcal {SD}_{n}^{UD}(k)$ can be written uniquely as a word $L=U^{x_1}D^{y_1}U^{x_2}D^{y_2}\cdots U^{x_k}D^{y_k}$, such that
\begin{equation*}
\begin{cases}
x_{1}+x_{2}+\cdots+ x_{k}=n, & x_1, x_2, \cdots, x_k\geq 1, \\
y_{1}+y_{2}+\cdots+y_{k}=n,  & y_1, y_2, \cdots, y_k\geq 1.
\end{cases}
\end{equation*}
The number of solutions for the $x_i$'s and for the $y_i$'s both equal to ${{n-k+k-1}\choose{k-1}}={{n-1}\choose{k-1}}$. Hence equation \eqref{udpeak} is proved.

Each $L^{\prime}\in\mathcal {SD}_{n}^{UU}(k)$ can be written uniquely as a word $L^{\prime}=U^{x_1}D^{y_1}U^{x_2}D^{y_2}\cdots U^{x_k}D^{y_k}U^{x_{k+1}}$,
such that
\begin{equation*}
\begin{cases}
x_{1}+x_{2}+\cdots+ x_{k}+x_{k+1}=n, & x_1, x_2, \cdots, x_{k+1}\geq 1 \\
y_{1}+y_{2}+\cdots+y_{k}=n,  & y_1, y_2, \cdots, y_k\geq 1
\end{cases}
\end{equation*}
There are ${{n-k+k+1-1}\choose{k}}={{n}\choose{k}}$ solutions for the $x_i$'s and
${{n-1}\choose{k-1}}$ solutions for the $y_i$'s. Hence equation \eqref{uupeak} is proved.

From \eqref{udpeak} and \eqref{uupeak} we have that the number of super Dyck paths with $k$ peaks of
order $n$ is
\begin{equation*}
{{n-1}\choose{k-1}}^{2}+{{n-1}\choose{k}}^{2}+2{{n-1}\choose{k-1}}{{n-1}\choose{k}}
={{n}\choose{k}}^{2}.
\end{equation*}
\qed

\begin{coro}
The number of Dyck paths of order $n$ with $k$ peaks is:
\[N(n,k)=\frac{1}{n}{{n}\choose{k}}{{n}\choose{k-1}}.\]
\end{coro}
\pf
From theorem \ref{th:main} we know that each Dyck path of
order $n$ with $k$ peaks is mapped to $k$ super Dyck paths, and each of the
$k$ super Dyck paths is either in $\mathcal{SD}_{n}^{UU}(k-1)$ or in $\mathcal{SD}
_{n}^{UD}(k)$. Therefore we have $kN(n,k)=\#\mathcal {S}\mathcal
{D}_{n}^{UU}(k-1)+\#\mathcal {S}\mathcal {D}_{n}^{UD}(k).$ From
Proposition \ref{yinli} we can conclude that \[
N(n,k)=\frac{1}{k}\left({{n-1}\choose{k-1}}^{2}
+{{n-1}\choose{k-2}}{{n-1}\choose{k-1}}\right)=\frac{1}{n}{{n}\choose{k}}{{n}\choose{k-1}}.
\]
\qed

Bijective proof of this result can also be found in \cite[Exercise 6.36(a)]{stanleyec2}.

\section{Humps in Schr\"{o}der paths}

In this section we count the number of humps in a third kind of lattice paths: Schr\"{o}der paths. A
\emph{Schr\"{o}der path} of order $n$ is a lattice path in
$\mathbb{Z}\times\mathbb{Z}$, from $(0,0)$ to $(n,n)$, using
up-steps $(0,1)$, down-steps $(1,0)$ and flat-steps $(1,1)$ (denoted
by $U$, $D$, $F$, respectively) and never going below the line $y=x$.
Note that Schr\"{o}der paths are different from rotating Motzkin paths 45 degrees counterclockwise, since the $F$ steps in these two kinds of paths are different. However, the bijection $\Phi$ still works when counting humps in Schr\"{o}der paths. Let $ss_n$ denote the number of super Schr\"{o}der paths of order $n$, and $hs_n$ denote the number of humps in all Schr\"{o}der paths of order $n$. We have the following result.

\begin{coro}\label{th:psn} For all $n \geq 0$, we have
\begin{equation}\label{sshs}
ss_{n}=2hs_{n}+1,
\end{equation}
and
\begin{equation}\label{hs}
hs_{n}=\frac{1}{2}\left(\sum_{k=0}^{n}{{n+k}\choose{2k}}{{2k}\choose{k}}-1\right).
\end{equation}
\end{coro}

\pf Apply the bijection $\Phi$ to Schr\"{o}der paths we immediately get \eqref{sshs}. Next we will count $ss_{n}$. Let $L$ be a super Schr\"{o}der path of order $n$ with $k$ humps, then there are $k$ $U$ steps,
$k$ $D$ steps, and $n-k$ $F$ steps in $L$. We can first choose a super Dyck path of order $k$ and then ``insert" $n-k$ $F$ steps to get $L$. There are ${2k \choose k}$ ways to choose a super Dyck paths, and  ${{n-k+2k+1-1}\choose{2k}}={{n+k}\choose{2k}}$ ways for
the insertion. Therefore we have
\[
ss_{n}=\sum_{k=0}^{n}{{n+k}\choose{2k}}{{2k}\choose{k}}.
\]
From the above formula and \eqref{sshs} we get \eqref{hs}.
\qed

The above proof inspired us to get the following identity, which is listed as an exercise in \cite[Exercise 3(g) of Chapter 1]{stanleynew}.

\begin{coro}For all $n \geq 0$, we have
\begin{equation}\label{schid}
\sum_{k=0}^{n} {{n}\choose{k}}^{2}(m+1)^{k}
=\sum_{k=0}^{n}{{n+k}\choose{2k}}{{2k}\choose{k}}m^{n-k}.
\end{equation}
\end{coro}
\pf We will first prove \eqref{schid} $m=1$. From the proof of Corollary \ref{th:psn} we
know that the right hand side of \eqref{schid} is the number of super
Schr\"{o}der paths of order $n$ when $m=1$ . Now we count $ss_{n}$ with a different
method to obtain the left hand. Let $L$ be a super Dyck path of order $n$ with $k$ peaks, for each peak of $L$, we can either keep it invariant or change it into a $F$ step to we get two super
Schr\"{o}der paths. Hence each $L$ is mapped to $2^{k}$ super Schr\"{o}der
paths, thus the left hand side of \eqref{schid} when $m=1$ also equals $ss_n$. Therefore we proved \eqref{schid} for $m=1$.

For general $m$ we count the number of super Schr\"{o}der paths in which the $F$ steps are $m$-colored. Now every super Dyck path with $k$ peaks is mapped to $(m+1)^k$ colored super Schr\"{o}der paths. So the total
number of such path is $\sum_{k=0}^{n}
{{n}\choose{k}}^{2}(m+1)^{k}$. On the other hand, from the proof of
Theorem \ref{th:psn} we know that the right hand side of \eqref{schid} also counts the number of such paths, hence we proved \eqref{schid}.
\qed

\vskip 2mm \noindent{\bf Acknowledgments.} This work is partially supported by the National Science Foundation of China under Grant No. 10801053, Shanghai Rising-Star Program (No. 10QA1401900), and the Fundamental Research Funds for the Central Universities.


\begin{thebibliography}{99}
\small \setlength{\itemsep}{-.8mm}

\bibitem{Deutsch} E. Deutsch, {\em Dyck path enumeration}, Discrete Math. 204 (1999), no. 1-3, 16--202.

\bibitem{Mansour1} T. Mansour, {\em Counting Peaks at Height $k$ in a Dyck Path}, Journal of Integer Sequences, Vol. 5 (2002), Article 02.1.1.

\bibitem{Mansour2} T. Mansour, {\em Statistics on Dyck Paths}, Journal of Integer Sequences, Vol. 9:1 (2006), Article 06.1.5.

\bibitem{Narayana} T.V. Narayana, {\em A partial order and its applications to probability}, Sankhya 21 (1959) 91--98.

\bibitem{PeartWoan} P. Peart and W. Woan, {\em Dyck Paths With No Peaks At Height $k$}, Journal of Integer Sequences, Vol. 4 (2001), Article 01.1.3.

\bibitem{A=B} M. Petkovesk, H. S. Wilf and D. Zeilberger, {em A=B}, AK Peters Ltd. (1996)

\bibitem{HumpsRegev} A. Regev, {\em Humps for Dyck and for Motzkin paths}, arXiv: 1002. 4504 v1
[math. CO] 24 Feb 2010.

\bibitem{sequence} N. J. A. Sloane, {\em Online Encyclopedia of Integer Sequence}, published electronically at http: //~oeis.~org

\bibitem{stanleynew}
R. P. Stanley, Enumerative Combinatorics (Volume 1 second
edition), http://www-math.mit.edu/$\thicksim$rstan/ec/ec1.pdf, 2011.

\bibitem{stanleyec1} R. P. Stanley, {\em Enumerative Combinatorics}, vol.~1, Cambridge Studies in Advanced Mathematics, vol. 49, Cambridge University Press, Cambridge, 1997.

\bibitem{stanleyec2} R. P. Stanley, {\em Enumerative Combinatorics}, vol.~2, Cambridge Studies in Advanced Mathematics, vol. 62, Cambridge University Press, Cambridge, 1999.

\bibitem{Zeilberger} D. Zeilberger,{\em The method of creative telescoping}, J. Symbolic Computation, 11 (1991) 195--204.

\end{thebibliography}
\end{document}